\documentclass[preprint,1p,times]{elsarticle}
\usepackage{amsmath, amssymb, amsthm}
\usepackage{graphicx}
\usepackage{color,soul}
\usepackage{float}
\usepackage{lineno}
\usepackage{pgf}
\usepackage{epstopdf}

\journal{ }
\begin{document}

\begin{frontmatter}
	
\title{Study of a Discretized Fractional-order eco-epidemiological model with prey infection}
\author[rmv]{Shuvojit Mondal\fnref{fn1}}
%\ead{shuvojitmondal91@gmail.com}
\author[chi]{Xianbing Cao\fnref{fn2}}
%\ead{}
\author[cmbe]{Nandadulal Bairagi\corref{cor1}}
\ead{nbairagi.math@jadavpuruniversity.in}
\cortext[cor1]{Corresponding author}
\fntext[cor1]{A part of this work was done at BTBU, Beijing, China when NB visited this university in February-March, 2019.}
\address[rmv]{Department of Mathematics, Rabindra Mahavidyalaya\\ Hooghly-712401, India.}
\address[chi]{Beijing Technology and Business University\\ Department of Mathematics, School of Science\\ Beijing, China.}
\address[cmbe]{Centre for Mathematical Biology and Ecology \\ Department of Mathematics, Jadavpur University\\ Kolkata-700032, India.}

\begin{abstract}
In this paper, an attempt is made to understand the dynamics of a three-dimensional discrete fractional-order eco-epidemiological model with Holling type II functional response. We first discretize a fractional-order predator-prey-parasite system with piecewise constant arguments and then explore the system dynamics. Analytical conditions for the local stability of different fixed points have been determined using the Jury criterion. Several examples are given to substantiate the analytical results. Our analysis shows that stability of the discrete fractional order system strongly depends on the step-size and the fractional order. More specifically, the critical value of the step-size, where the switching of stability occurs, decreases as the order of the fractional derivative decreases. Simulation results explore that the discrete fractional-order system may also exhibit complex dynamics, like chaos, for higher step-size. 
\end{abstract}


\end{frontmatter}

 \section{Introduction}\label{sec:1}

     The idea of fractional calculus has been known since the notion of fractional derivative was firstly introduced by Leibniz in a letter dated 30th September, 1695, where half-order derivative was mentioned. Since then it was subsequently developed by Euler, Lagrange, Laplace, Fourier, Abel, Liouville, Riemann, Heaviside, Caputo along with many others \cite{SET93}. Generally, fractional calculus deals with the study of fractional order integral and derivative operators over real or complex domains and their applications \cite{Rihan15}. Specially it is a generalization of classical differential and integral calculus of integer order to arbitrary order. In recent times, fractional order derivatives and fractional order differential equations have been applied in several fields of science and engineering \cite{MET11,MO15,GET10,AET12,D11,M06,LZ18,M13}, however, initially it was treated as a topic of interest of pure mathematicians only \cite{TET15}. The reason is two-fold.  First, fractional derivatives have an additional degree of freedom over its integer order counterpart due to the additional parameter that represents its order and more suitable for those systems having higher order dynamics and complex nonlinear phenomena \cite{TB84,SET11}. Secondly and more importantly, fractional order derivatives not only depend on the local conditions but also on the history of the function \cite{BET17} and, therefore, fractional derivatives have become an efficient tool for those systems, where the consideration of memory or hereditary properties of the function is essential to represent the system, e.g., in the case of biological systems. In the last two decades, fractional order calculus has found its many applications in biological sciences \cite{Ahmed07,RET13,CY14,MET19,MET17,HET15, LET16,V15,MET17a}.\\

Discrete system of a continuous time model has significant applications in simulations and computations. Particularly, the population model can be better described in discrete time and it is especially applicable in populations having non-overlapping generations, e.g., insect populations. Such discrete models are known to show richer dynamics compare to its continuous counterpart \cite{Hattaf15, Sekiguchi10, Chen13, Hu12, Biswas17,GET20}. Thus, both from biological and dynamic point of views, study of discrete system is important. Though the number of discrete population models is numerous, the study of discrete fractional-order population models is rarely observed in the literature. Recently, a transformed method was employed to obtain solutions of a family of finite fractional difference equations in \cite{AE07}. The global and local existence of solutions for the fractional difference equation are given in \cite{CET1107}. Numerical solution was adopted to synchronize chaos in a discrete fractional logistic map \cite{WB14}. For an ecological system, Elsadany and Matouk \cite{Matouk15} recently studied a two-dimensional fractional-order Lotka-Volterra predator prey model with its discretization and showed complex dynamics. The existence and uniqueness of solutions in two-dimensional discrete fractional Lotka-Volterra model was addressed in \cite{Alzabut18}. Selvam and Janagaraj \cite{SJ18} studied the dynamics of a two-dimensional discrete fractional order predator-prey system with type II functional response. A two-dimensional fractional-order predator-prey model with saturated harvesting was studied in \cite{SET18}. Local stability analysis and bifurcation analysis have been done in a simple fractional-order SIR model by Salman \cite{Salman17}. Similar analysis was done for a fractional-order discrete SI epidemic model in \cite{AET18}. Discretization has also been done in other fractional-order biological and physical systems \cite{Khan19,Khan20,Mondal20,Khan2020,Agarwal13,Sayed13,Selvam18}.
These relatively simple models show that the fractional-order plays a critical role in breaking the stability. By considering a simple two-dimensional epidemic model, Raheem and Salman \cite{Raheem14} have shown that the dynamic complexity depends also on the step-size of the discrete system. Analytical complexity, however, multiplies if the dimension of the system increases and/or strong nonliearity is there. In this paper, we construct a discrete version of a highly nonlinear three-dimensional fractional-order  predator-prey-parasite (PPP) model and reveal its dynamics. So far our knowledge goes, no study has been done on a discrete fractional-order PPP system. We unveil the dynamic properties of this highly nonlinear model both from analytical and numerical points of view. We show analytically that the stability of different fixed points depends on both the fractional order and step-size. It also shows complex dynamics like chaos as shown by the continuous time integer order system for some parameter values.

For this study, we assume the following continuous-time PPP model studied by Chattopadhyay and Bairagi \cite{ChattoBairagi01}: 

\begin{eqnarray}\label{Eco-epidemiological integer order model}
\frac{dX}{dt} & = & rX\bigg(1-\frac{X+Y}{K}\bigg) - \lambda YX,\nonumber \\
\frac{dY}{dt}& = & \lambda YX-\frac{mYZ}{a+Y} - \mu Y, \\
\frac{dZ}{dt}& = & \frac{\theta YZ}{a+Y}-dZ. \nonumber
\end{eqnarray}
It says that a prey species is divided into two sub classes: susceptible prey $(X)$ and infected prey $(Y)$ due to some micro-parasite infection. Susceptible prey can give birth and grow logistically to its maximum value $K$ supported by the environment with intrinsic growth rate $r$. The infection spreads by contact following mass action law with $\lambda$ as its force of infection. Infected prey is unable to reproduce but share resources with susceptible prey. Infected prey dies out at a rate $\mu$, which considers both the natural and infection-related deaths. Predator population $Z$ consumes infected preys only as they are weakened due to infection and can not escape predation. Here $m$ is the maximum prey attack rate of predator and $a$ is the half-saturation constant. The parameter $\theta$ ($0<\theta<1$) measures the reproductive gain of predator and $d$ measures its natural mortality. All parameters are assumed to be positive from biological point of view. Readers are referred to \cite{ChattoBairagi01} for more discussion about the model.

%\newpage       
\begin{center}{Table 1.1. Parameter descriptions and their units}
\end{center}

\hspace{-0.8cm}
{\small
	\begin{tabular}
		{l l l} \hline 
		Parameters & Description & Units  \\ \hline 
		
		$X$                 & Density of susceptible prey &  Number per unit designated area \\
		$Y$                 & Density of infected prey &  Number per unit designated area \\
		$Z$                 & Density of predator &  Number per unit designated area \\
		$r$                 & Intrinsic growth rate of prey&  Per day \\ 
		$K$                 & Carrying capacity & Number per unit designated area \\ 
		$a$                 & Half-saturation constant & Number per unit designated area \\
		$\lambda$                 & Force of infection & Per day \\ 
		$m$                 & Search rate for infected prey & Per day\\	
		$\mu$               & Total (natural + disease-induced) & Per day\\
		& death rate of infected prey & \\	
		$\theta$            & Reproductive gain & Per day\\	
		$d$                 & Death rate of predator & Per day\\
		$\alpha$            & Fractional order  & $0 < \alpha < 1$ \\
		\hline
		
\end{tabular}}
\vspace{.2in}      

Considering the fractional derivative in caputo sense, Mondal et al. \cite{MET17a} extended the work of Chattopadhyay and Bairagi \cite{ChattoBairagi01} and analyzed the following fractional order PPP model (\ref{Eco-epidemiological integer order model}):
\begin{eqnarray}\label{Eco-epidemiological fractional order model}
^{c}_{0} D^{\alpha}_{t}X & = & rX\bigg(1-\frac{X+Y}{K}\bigg) - \lambda YX,~X(0) > 0, \nonumber \\
^{c}_{0} D^{\alpha}_{t}Y & = & \lambda YX-\frac{mYZ}{a+Y} - \mu Y,~Y(0) > 0, \\
^{c}_{0} D^{\alpha}_{t}Z & = & \frac{\theta YZ}{a+Y}-dZ,~Z(0) > 0, \nonumber
\end{eqnarray}
where $^{c}_{0} D^{\alpha}_{t}$ is the Caputo fractional derivative with fractional order $\alpha$ $(0< \alpha < 1)$. It is worth mentioning that the systems (\ref{Eco-epidemiological integer order model}) and (\ref{Eco-epidemiological fractional order model}) become identical for $\alpha=1$. They have shown that the solutions of system (\ref{Eco-epidemiological fractional order model}) are positively invariant and uniformly bounded in $R^3_+$ under some restrictions. Existence and uniqueness of solution of system (\ref{Eco-epidemiological fractional order model}) have also been discussed. Local and global stability of all biologically feasible equilibrium points of system (\ref{Eco-epidemiological fractional order model}) were also proved. Here we discritize the model system (\ref{Eco-epidemiological fractional order model}) and explore its dynamics.

\section{Discretized fractional-order model and its analysis}
Following Elsadany and Matouk \cite{Matouk15}, discretization of model system (\ref{Eco-epidemiological fractional order model}) with piecewise constant arguments can be done in the following manner:
\begin{eqnarray}\nonumber
^{c}_{0} D^{\alpha}_{t}X & = & rX([t/s]s)\bigg(1-\frac{X([t/s]s)+Y([t/s]s)}{K}\bigg) - \lambda X([t/s]s)Y([t/s]s), \\
^{c}_{0} D^{\alpha}_{t}Y & =  & \lambda X([t/s]s) Y([t/s]s) - \frac{m Y([t/s]s) Z([t/s]s)}{(a + Y([t/s]s))} - \mu Y([t/s]s),\nonumber\\
^{c}_{0} D^{\alpha}_{t}Z & =  &  \frac{\theta Y([t/s]s) Z([t/s]s)}{(a + Y([t/s]s))} - d Z([t/s]s),\nonumber
\end{eqnarray}
with initial condition $X(0) = X_0 > 0$, $Y(0) = Y_0> 0, Z(0) = Z_0 >0$ and $s (>0)$ is the step-size. Let $t\in [0,s)$, so that $t/s \in [0,1)$. In this case, we have
\begin{eqnarray}\nonumber
^{c}_{0} D^{\alpha}_{t}X & = & X_0\bigg(r(1-\frac{X_0+Y_0}{K}) - \lambda Y_0 \bigg), \\
^{c}_{0} D^{\alpha}_{t}Y & =  & Y_0\bigg(\lambda X_0 - \frac{m Z_0}{(a + Y_0)} - \mu\bigg),\nonumber\\
^{c}_{0} D^{\alpha}_{t}Z & =  &  Z_0\bigg(\frac{\theta Y_0}{(a + Y_0)} - d\bigg),\nonumber
\end{eqnarray}
and the solution of this fractional differential equation can be written as 
\begin{eqnarray}\nonumber
X_1(t) & = & X_0 + J^{\alpha}_{0}\bigg(X_0\bigg(r(1-\frac{X_0 + Y_0}{K}) - \lambda Y_0\bigg)\bigg) \nonumber \\
&  = & X_0 + \frac{t^\alpha}{\alpha\Gamma (\alpha)} \bigg(X_0\bigg(r(1-\frac{X_0 + Y_0}{K}) - \lambda Y_0\bigg)\bigg), \nonumber \\
Y_1(t) & = & Y_0 + J^{\alpha}_{0}\bigg(Y_0\bigg(\lambda X_0 - \frac{m Z_0}{(a + Y_0)} - \mu\bigg)\bigg) \nonumber \\
& = & Y_0 + \frac{t^\alpha}{\alpha \Gamma (\alpha)}\bigg(Y_0\bigg(\lambda X_0 - \frac{m Z_0}{(a + Y_0)} - \mu\bigg)\bigg),\nonumber\\
Z_1(t) & = & Z_0 + J^{\alpha}_{0}\bigg(Z_0\bigg(\frac{\theta Y_0}{(a + Y_0)} - d\bigg)\bigg) \nonumber \\
& = & Z_0 + \frac{t^\alpha}{\alpha \Gamma (\alpha)}\bigg(Z_0\bigg(\frac{\theta Y_0}{(a + Y_0)} - d\bigg)\bigg),\nonumber
\end{eqnarray}
where $J^\alpha_{0} = \frac{1}{\Gamma(\alpha)}\int^{t}_{0} (t- \tau)^{\alpha-1} d\tau, \alpha>0$.

\noindent In the second step, we assume $t\in [s,2s)$ so that $t/s \in [1,2)$ and similarly obtain
\begin{eqnarray}\nonumber
^{c}_{0} D^{\alpha}_{t}X & = & X_1(s)\bigg(r(1-\frac{X_1(s)+Y_1(s)}{K}) - \lambda Y_1(s) \bigg), \\
^{c}_{0} D^{\alpha}_{t}Y & =  & Y_1(s)\bigg(\lambda X_1(s) - \frac{m Z_1(s)}{(a + Y_1(s))} - \mu\bigg),\nonumber\\
^{c}_{0} D^{\alpha}_{t}Z & =  &  Z_1(s)\bigg(\frac{\theta Y_1(s)}{(a + Y_1(s))} - d\bigg). \nonumber
\end{eqnarray}

\noindent The solution of this equation reads
\begin{eqnarray}\nonumber
X_2(t) & = & X_1(s) + J^{\alpha}_{s}\bigg(X_1(s)\bigg(r(1-\frac{X_1(s)+Y_1(s)}{K}) - \lambda Y_1(s)\bigg)\bigg) \nonumber\\
& = & X_1(s) + \frac{(t - s)^\alpha}{\alpha\Gamma (\alpha)}\bigg(X_1(s)\bigg(r(1-\frac{X_1(s)+Y_1(s)}{K}) - \lambda Y_1(s)\bigg)\bigg), \nonumber
\end{eqnarray}
\begin{eqnarray}\nonumber
Y_2(t) & = & Y_1(s) + J^{\alpha}_{s}\bigg(Y_1(s)\bigg(\lambda X_1(s) - \frac{m Z_1(s)}{(a + Y_1(s))} - \mu\bigg)\bigg) \nonumber\\
& = & Y_1(s) + \frac{(t - s)^\alpha}{\alpha\Gamma (\alpha)}\bigg(Y_1(s)\bigg(\lambda X_1(s) - \frac{m Z_1(s)}{(a + Y_1(s))} - \mu\bigg)\bigg),\nonumber\\
Z_2(t) & = & Z_1(s) + J^{\alpha}_{s}\bigg(Z_1(s)\bigg(\frac{\theta Y_1(s)}{(a + Y_1(s))} - d\bigg)\bigg) \nonumber\\
& = & Z_1(s) + \frac{(t - s)^\alpha}{\alpha\Gamma (\alpha)}\bigg(Z_1(s)\bigg(\frac{\theta Y_1(s)}{(a + Y_1(s))} - d\bigg)\bigg),\nonumber
\end{eqnarray}
where $J^\alpha_{s} = \frac{1}{\Gamma(\alpha)}\int^{t}_{s} (t- \tau)^{\alpha-1} d\tau, \alpha>0$.\\ Repeating the discretization process $n$ times, we have
\begin{eqnarray}\nonumber
X_{n+1}(t) & = & X_n(ns) + \frac{(t - ns)^\alpha}{\alpha\Gamma (\alpha)}\bigg(X_n(ns)\bigg(r(1-\frac{X_n(ns) + Y_n(ns)}{K}) - \lambda Y_n(ns)\bigg)\bigg), \nonumber \\
Y_{n+1}(t) & = & Y_n(ns) + \frac{(t - ns)^\alpha}{\alpha\Gamma (\alpha)}\bigg(Y_n(ns)\bigg(\lambda X_n(ns) - \frac{m Z_n(ns)}{(a + Y_n(ns))} - \mu\bigg)\bigg), \nonumber\\
Z_{n+1}(t) & = & Z_n(ns) + \frac{(t - ns)^\alpha}{\alpha\Gamma (\alpha)}\bigg(Z_n(ns)\bigg(\frac{\theta Y_n(ns)}{(a + Y_n(ns))} - d\bigg)\bigg), \nonumber
\end{eqnarray}
where $t\in [ns, (n+1)s)$.\\ Making $t\rightarrow (n+1)s$, the corresponding fractional discrete model of the continuous fractional model (\ref{Eco-epidemiological fractional order model}) is obtained as
\begin{eqnarray}\nonumber \label{Discrete model}
X_{n+1} & = & X_n + \frac{s^\alpha}{\alpha\Gamma (\alpha)}\bigg(X_n\bigg(r(1-\frac{X_n + Y_n}{K}) - \lambda Y_n\bigg)\bigg), \nonumber \\
Y_{n+1} & = & Y_n + \frac{s^\alpha}{\alpha\Gamma (\alpha)}\bigg(Y_n\bigg(\lambda X_n - \frac{m Z_n}{(a + X_n)} - \mu\bigg)\bigg), \nonumber\\
Z_{n+1} & = & Z_n + \frac{s^\alpha}{\alpha\Gamma (\alpha)}\bigg(Z_n\bigg(\frac{\theta Y_n}{(a + Y_n)} - d\bigg)\bigg).
\end{eqnarray}
For $\alpha = 1$, one obtains the so called Euler forward discrete model as a special case of this generalized discrete model.

\subsection{Existence and stability of fixed points}\label{subsec:2.1}

In the following, we investigate the dynamics of the discretized fractional order model (\ref{Discrete model}). At the fixed point, we have $X_{n+1} = X_{n}$, $Y_{n+1} = Y_{n}$ and $Z_{n+1} = Z_{n}$. One can then easily compute that (\ref{Discrete model}) has the same fixed points as in the integer and fractional order systems (\ref{Eco-epidemiological integer order model}) and (\ref{Eco-epidemiological fractional order model}). The system (\ref{Discrete model}) has four equilibrium points \cite{MET17a, ChattoBairagi01}, viz. $E_{0} = (0,0,0)$ as the trivial equilibrium, $E_{1} = (K,0,0)$ as the axial equilibrium, $E_{2} = (X_1, Y_1, 0)$ as the planar equilibrium, where $X_1 = \frac{\mu}{\lambda}$ and $Y_1 = \frac{r(\lambda K-\mu)}{\lambda(r + \lambda K)}$ and $E^{*} = (X^{*},Y^{*},Z^{*})$ as the interior equilibrium, where
\begin{eqnarray}\label{Equilibrium relation}
X^{*} =  K-\bigg(1 + \frac{\lambda K}{r}\bigg)Y^{*},~~ Y^{*} =  \frac{ad}{\theta -d},~~ Z^{*} = \frac{(a+Y^{*})(\lambda X^{*} - \mu)}{m}.
\end{eqnarray}
Note that the equilibria $E_{0}$ and $E_{1}$ always exist. The planar equilibrium point $E_{2}$ exists if $R_{0}>1$, where $R_{0} =\frac{\lambda K}{\mu}$. The interior equilibrium $E^{*}$ exists if $(i)$ $R_{0}>1$ and $(ii)$ $\theta > \theta_{1}$, where $\theta_1=d + \frac{\lambda ad(r + \lambda K)}{r(\lambda K - \mu)}$.\\

The Jacobian matrix of system (\ref{Discrete model}) at any arbitrary fixed point $(X, Y, Z)$ reads
\begin{equation}
{J(X, Y, Z)}=
\begin{pmatrix}
a_{11} & a_{12} & 0 \\ a_{21} & a_{22} & a_{23} \\ 0 & a_{32} & a_{33}
\end{pmatrix},
\label{eq:myeqn}
\end{equation}
where
\begin{eqnarray}\nonumber
a_{11} & =& 1 + \frac{s^\alpha}{\alpha\Gamma (\alpha)}\bigg(r(1-\frac{2X + Y}{K}) - \lambda Y\bigg),\nonumber\\
a_{12} & = & - \frac{s^\alpha}{\alpha\Gamma (\alpha)} \bigg(X(\lambda + \frac{r}{K})\bigg), \nonumber \\
a_{21} & = & \frac{s^\alpha}{\alpha\Gamma (\alpha)} \lambda Y,\nonumber \\
a_{22} & =& 1 + \frac{s^\alpha}{\alpha\Gamma (\alpha)} \bigg(\lambda X - \frac{mZ}{a+Y}-  \mu\bigg) + \frac{s^\alpha}{\alpha\Gamma (\alpha)} \frac{mYZ}{(a+Y)^2},\nonumber\\
a_{23} & =& -\frac{s^\alpha}{\alpha\Gamma (\alpha)} \frac{mY}{a+Y},\nonumber\\
a_{32} & =& \frac{s^\alpha}{\alpha\Gamma (\alpha)} \frac{a\theta Z}{(a+Y)^2},\nonumber\\
a_{33} & =& 1 + \frac{s^\alpha}{\alpha\Gamma (\alpha)} \bigg(\frac{\theta Y}{a+Y} - d\bigg).\nonumber
\end{eqnarray}
Let $\xi_1$, $\xi_2$ and $\xi_3$ be the eigenvalues of the Jacobian matrix (\ref{eq:myeqn}). Then we have the following definition and lemma.

\noindent\textbf{Definition 1} \cite{Salman17, Selvam18} A fixed point $(X, Y, Z)$ of system (\ref{Discrete model}) is called a sink which is locally asymptotically stable if $\mid\xi_i\mid < 1$ $(i=1,2,3)$ for all $i$. It is called a source which is unstable if $\mid\xi_i\mid >1$ $(i=1,2,3)$ for all $i$. It is called a saddle if any two eigenvalues follow opposite inequality, and a non-hyperbolic fixed point if $\mid\xi_i\mid = 1$ $(i=1,2,3)$ for at least one $i$.

\noindent \textbf{Lemma 1} \cite{Elaydi99} Suppose the characteristic polynomial $F(\xi)$ of a second order Jacobian matrix is given by $F(\xi) = \xi^2 - a_1 \xi + a_2$. Then the solutions $\xi_i,~ i = 1, 2,$ of $F(\xi) = 0$ satisfy  $\mid\xi_i\mid < 1$ for $i = 1, 2,$ if the following conditions hold: $$(i)~a_2 < 1,~ (ii)~ 1 + a_2 > \mid a_1 \mid.$$

\noindent \textbf{Lemma 2} \cite{Elaydi99} Suppose the characteristic polynomial $p(\xi)$ of Jacobian matrix (\ref{eq:myeqn}) is given by $p(\xi) = \xi^3 + a_1 \xi^2 + a_2 \xi + a_3$. Then the solutions $\xi_i,~ i = 1, 2, 3,$ of $p(\xi) = 0$ satisfy  $\mid\xi_i\mid < 1$ for $i = 1, 2, 3$ if the following conditions hold: \\$(i)~p(1) = 1 + a_1 + a_2 + a_3 > 0,~ (ii)~ (-1)^3p(-1) = 1 - a_1 + a_2 - a_3 > 0,~ \mbox{and} ~(iii)~ 1 - (a_3)^2 > \mid a_2 -a_3a_1 \mid.$

\subsection{Main result}
We now prove our main theorem in relation to the stability of the discrete fractional order system (\ref{Discrete model}).\\

\noindent \textbf{Theorem 1} (a) The fixed point $E_0$ is always unstable for any feasible values of  $\alpha$ and $s$. \\
	%It will be a saddle point if $s_1 < s < s_2$ and a source if $s > max \{s_1, s_2\}$ for any $\alpha \in (0, 1]$. If either $s = s_1$ or $s= s_2$, then $E_0$ is nonhyperbolic, where $s_1 = \sqrt[\alpha]{\frac{2\alpha\Gamma(\alpha)}{\mu}}, ~ s_2 = \sqrt[\alpha]{\frac{2\alpha\Gamma(\alpha)}{d}}$.  \\
	
	\noindent (b) The fixed point $E_1$ is locally asymptotically stable for $0<\alpha\leq 1$ if $R_0 = \frac{\lambda K}{\mu} <1$, $s < min \{s_2 , s_3, s_4\}$, where $s_2 = \sqrt[\alpha]{\frac{2\alpha\Gamma(\alpha)}{d}}$, $s_3 = \sqrt[\alpha]{\frac{2\alpha\Gamma(\alpha)}{r}}$ and $s_4 = \sqrt[\alpha]{\frac{2\alpha\Gamma(\alpha)}{\mu - \lambda K}}$.  When $R_0 > 1$,  the fixed point $E_1$ is always unstable for $0<\alpha\leq 1$ for any step size $s$. \\
	%It is a saddle point if $R_0  >1$, $s_3 < s < s_2$~; or $R_0>1$, $s_2 < s < s_3$; and a source if $R_0 > 1$, $s > max\{s_2, s_3\}$ for any $\alpha \in (0, 1]$. \\
	
	\noindent (c) The fixed point $E_2$ is locally asymptotically stable for $0<\alpha\leq 1$ if $R_0 = \frac{\lambda K}{\mu} >1$, $ s_5< s < min \{s_6 , s_7\}$ with $d> d_1$, where $s_5 = \sqrt[\alpha]{\frac{\alpha\Gamma(\alpha)}{\lambda K -\mu}}$, $s_6 = \sqrt[\alpha]{\frac{2\alpha\Gamma(\alpha)}{d-d_1}}$, $s_7 = \sqrt[2\alpha]{\frac{\lambda K (\alpha\Gamma(\alpha))^2}{\mu r (\lambda K - \mu)}}$ and $d_1 = \frac{\theta r (\lambda K - \mu)}{a\lambda(\lambda K + r)+ r(\lambda K-\mu)}$. When $d < d_1$,  the fixed point $E_2$ is always unstable for $0<\alpha\leq 1$.\\
	
	\noindent (d) The fixed point $E^*$ is locally asymptotically stable for $0<\alpha\leq 1$ if $R_0 > 1$, $0 < s < min\{s_{8}, s_{9}\}$ with $\theta > \theta_{1}$, where  $s_8 = \sqrt[\alpha]{\frac{2K\alpha\Gamma(\alpha)}{rX^*}}$, $\theta_1=d + \frac{\lambda ad(r + \lambda K)}{r(\lambda K - \mu)}$ and  $s_{9}$ is such that \eqref{Jury} holds.

\noindent \textbf{Proof}
	At the fixed point $E_0$, the Jacobian matrix $J(E_{0})$ can be obtained from (\ref{eq:myeqn}) and it is given by
	\[
	J(E_{0}) = \begin{pmatrix}
	1+\frac{s^\alpha}{\alpha\Gamma (\alpha)}r & 0 & 0 \\ 0 & 1-\frac{s^\alpha}{\alpha\Gamma (\alpha)}\mu & 0 \\ 0 & 0 & 1-\frac{s^\alpha}{\alpha\Gamma (\alpha)}d \end{pmatrix}.\]
	The eigenvalues are $\xi_1 = 1+\frac{s^\alpha}{\alpha\Gamma (\alpha)}r$, $\xi_2 =  1-\frac{s^\alpha}{\alpha\Gamma (\alpha)} \mu$ and $ \xi_3 = 1 - \frac{s^\alpha}{\alpha\Gamma (\alpha)}d$. Since $|\xi_1| > 1$, $E_0$ is always unstable for $0<\alpha \leq 1$ and any $s>0$. \\
	%In fact, it is a saddle point if one of the following conditions hold: $(i)~ \sqrt[\alpha]{\frac{2\alpha\Gamma(\alpha)}{\mu}} < s < \sqrt[\alpha]{\frac{2\alpha\Gamma(\alpha)}{d}}$ for which $|\xi_2 | > 1$, $|\xi_3 | < 1$;  $(ii)~\sqrt[\alpha]{\frac{2\alpha\Gamma(\alpha)}{d}} < s < \sqrt[\alpha]{\frac{2\alpha\Gamma(\alpha)}{\mu}}$ for which  $|\xi_2 | < 1$, $|\xi_3 | > 1$ and $(iii)~ 0 < s < min\biggl\{\sqrt[\alpha]{\frac{2\alpha\Gamma(\alpha)}{\mu}}, \sqrt[\alpha]{\frac{2\alpha\Gamma(\alpha)}{d}}\biggr\}$ for which $|\xi_2 | < 1$, $|\xi_3 | < 1$. Again it is a source if $s > max\{\sqrt[\alpha]{\frac{2\alpha\Gamma(\alpha)}{\mu}}, \sqrt[\alpha]{\frac{2\alpha\Gamma(\alpha)}{d}}\}$ for which $|\xi_{2,3}| > 1$. Again, it becomes nonhyperbolic if either $s = \sqrt[\alpha]{\frac{2\alpha\Gamma(\alpha)}{\mu}}$ or $s = \sqrt[\alpha]{\frac{2\alpha\Gamma(\alpha)}{d}}$ for any $\alpha \in (0,1]$. \\
	
	\noindent The Jacobian matrix $J(E_{1})$ at the fixed point $E_1$ is evaluated as
	\[
	J(E_{1}) = \begin{pmatrix}
	1-\frac{s^\alpha}{\alpha\Gamma (\alpha)}r & \frac{s^\alpha}{\alpha\Gamma (\alpha)}(\lambda K +r) & 0 \\ 0 & 1 + \frac{s^\alpha}{\alpha\Gamma (\alpha)}(\lambda K -\mu) & 0 \\ 0 & 0 & 1-\frac{s^\alpha}{\alpha\Gamma (\alpha)}d \end{pmatrix}.\]
	The eigenvalues of $J(E_{1})$ are 
	\begin{equation}\nonumber
	\xi_1 = 1-\frac{s^\alpha}{\alpha\Gamma (\alpha)}r,~ \xi_2 = 1 + \frac{s^\alpha}{\alpha\Gamma (\alpha)}(\lambda K -\mu),~ \xi_3 = 1 - \frac{s^\alpha}{\alpha\Gamma (\alpha)}d. \nonumber
	\end{equation}
	
	As $0<\alpha \leq 1$,  $|\xi_{1,2,3} | < 1$ hold if $$s < min \bigg\{\sqrt[\alpha]{\frac{2\alpha\Gamma(\alpha)}{r}}, \sqrt[\alpha]{\frac{2\alpha\Gamma(\alpha)}{d}}, \sqrt[\alpha]{\frac{2\alpha\Gamma(\alpha)}{\mu - \lambda K}}\bigg\}, \mu>\lambda k.$$ 
	Therefore, $E_1$ is locally asymptotically stable for $0<\alpha\leq 1$ if $R_0 = \frac{\lambda K}{\mu} <1$, $s < min \{s_2 , s_3, s_4\}$, where $s_2 = \sqrt[\alpha]{\frac{2\alpha\Gamma(\alpha)}{d}}$, $s_3 = \sqrt[\alpha]{\frac{2\alpha\Gamma(\alpha)}{r}}$ and $s_4 = \sqrt[\alpha]{\frac{2\alpha\Gamma(\alpha)}{\mu - \lambda K}}$. However, if $R_0 > 1$ then $|\xi_2| > 1$. So $E_{1}$ is always unstable for all $0<\alpha\leq 1$ and for any step size $s$.\\
	% Note that $|\xi_1| > 1$ if $s > s_2$ and $|\xi_3| > 1$ if $s > s_3$. Thus, $E_1$ will be a source if $R_0 > 1$, $s > max\{s_2, s_3\}$ for any $\alpha \in (0, 1]$. The fixed point $E_1$ will be a saddle point if \textcolor{red}{ one of the conditions $(i)~ R_0  >1,~ s_2 < s < s_3$~; or  $(ii)~ R_0>1,~ s_3 < s < s_2$~;  or  $(iii)~ R_0>1,~ 0 < s < min\{s_2, s_3\}$ hold.}\\
	
	\noindent At the fixed point $E_2$, the Jacobian matrix $J(E_2)$ takes the form 
	\[
	J(E_{2}) = \begin{pmatrix}
	1-\frac{s^\alpha}{\alpha\Gamma (\alpha)}\frac{r\mu}{\lambda K} & -\frac{s^\alpha}{\alpha\Gamma (\alpha)}\frac{\mu(r + \lambda K)}{\lambda K} & 0 \\ \frac{r(\lambda K - \mu)}{(\lambda K + r)}\frac{s^\alpha}{\alpha\Gamma (\alpha)} & 0 & -\frac{m r (\lambda K - \mu)}{a\lambda(\lambda K + r)+ r(\lambda K-\mu)}\frac{s^\alpha}{\alpha\Gamma (\alpha)} \\ 0 & 0 & 1 + \frac{s^\alpha}{\alpha\Gamma (\alpha)}\bigg(\frac{\theta r (\lambda K - \mu)}{a\lambda(\lambda K + r)+ r(\lambda K-\mu)} - d\bigg) \end{pmatrix}.\]
	
	\noindent After some mathematical manipulations, one can have the following characteristic equation
	\begin{equation}\label{chareqn}
	\bigg(1 + \frac{s^\alpha}{\alpha\Gamma (\alpha)}\bigg(\frac{\theta r (\lambda K - \mu)}{a\lambda(\lambda K + r)+ r(\lambda K-\mu)} - d\bigg) - \xi\bigg)( \xi^{2} + A\xi + B)=0,
	\end{equation}
	where $A = \frac{s^\alpha}{\alpha\Gamma (\alpha)}\frac{r\mu}{\lambda K} - 1$ and $B = \bigg(\frac{s^\alpha}{\alpha\Gamma (\alpha)}\bigg)^2 \frac{r\mu (\lambda K- \mu)}{\lambda K}$.
	Therefore, one eigenvalue is $\xi_{1} = 1 + \frac{s^\alpha}{\alpha\Gamma (\alpha)} (d_{1} - d)$, where $d_{1} = \frac{\theta r (\lambda K - \mu)}{a\lambda(\lambda K + r)+ r(\lambda K-\mu)}$ and the other two eigenvalues $\xi_{2,3}$ are the roots of $\xi^{2} + A\xi + B = 0$. Following Lemma $1$, we observe $\mid \xi_{2,3} \mid < 1$ if  $ \sqrt[\alpha]{\frac{\alpha\Gamma(\alpha)}{\lambda K -\mu}} < s < \sqrt[2\alpha]{\frac{\lambda K (\alpha\Gamma(\alpha))^2}{\mu r (\lambda K - \mu)}}$ with $R_0 > 1$. For $\xi_1$, we consider the following two cases:\\
	
	\noindent\textbf{Case I:} If $d > d_{1}$ then $\xi_{1} = 1 - \frac{s^\alpha}{\alpha\Gamma (\alpha)} (d - d_1)$ and $\mid \xi_1 \mid < 1$ if $0 < s <\sqrt[\alpha]{\frac{2\alpha\Gamma(\alpha)}{d-d_1}}$. Therefore, the equilibrium $E_{2}$ is locally asymptotically stable for $0<\alpha\leq 1$ if $R_0 = \frac{\lambda K}{\mu} >1$, $ s_5< s < min \{s_6 , s_7\}$ with $d> d_1$, where $s_5 = \sqrt[\alpha]{\frac{\alpha\Gamma(\alpha)}{\lambda K -\mu}}$, $s_6 = \sqrt[\alpha]{\frac{2\alpha\Gamma(\alpha)}{d-d_1}}$, $s_7 = \sqrt[2\alpha]{\frac{\lambda K (\alpha\Gamma(\alpha))^2}{\mu r (\lambda K - \mu)}}$ and $d_1 = \frac{\theta r (\lambda K - \mu)}{a\lambda(\lambda K + r)+ r(\lambda K-\mu)}$.\\
	
	\noindent\textbf{Case II:} If $d < d_{1}$ then $\xi_{1} = 1 + \frac{s^\alpha}{\alpha\Gamma (\alpha)} (d_1 - d) > 1$. So $E_{2}$ is always unstable for all $0<\alpha\leq 1$ and for any step size $s$.\\

	\noindent At the interior fixed point $E^*$, the Jacobian matrix $J(E^*)$ is evaluated as
	\[
	{J(E^*)} = \begin{pmatrix}
	1-\frac{s^\alpha}{\alpha\Gamma (\alpha)}\frac{rX^*}{K} & -\frac{s^\alpha}{\alpha\Gamma (\alpha)}X^*(\frac{r}{K} + \lambda) & 0 \\ \frac{s^\alpha}{\alpha\Gamma (\alpha)} \lambda Y^*&  1 +\frac{s^\alpha}{\alpha\Gamma (\alpha)} \frac{mY^*Z^*}{(a+Y^*)^2} & -\frac{m Y^*}{a + Y^*}\frac{s^\alpha}{\alpha\Gamma (\alpha)} \\ 0 & \frac{s^\alpha}{\alpha\Gamma (\alpha)} \frac{\theta a Z^*}{(a+Y^*)^2} & 1 \end{pmatrix}.\]
	After some algebraic manipulation, the characteristic equation of $J(E^*)$ can be expressed as
	\begin{equation}\label{chareqn_2}
	p(\xi) = \xi^{3} + A_1\xi^2 + A_2\xi + A_3=0,
	\end{equation}
	where
	\begin{eqnarray} \nonumber
	A_1 & = & \frac{s^\alpha}{\alpha\Gamma (\alpha)}\frac{rX^*}{K} -  \frac{s^\alpha}{\alpha\Gamma (\alpha)} \frac{mY^*Z^*}{(a+Y^*)^2} - 3, \nonumber \\
	A_2 & = & 3 +  \frac{2s^\alpha}{\alpha\Gamma (\alpha)}\bigg(\frac{mY^*Z^*}{(a+Y^*)^2} - \frac{rX^*}{K}\bigg) + \bigg(\frac{s^\alpha}{\alpha\Gamma (\alpha)}\bigg)^2 \bigg(\frac{amdZ^*}{(a+Y^*)^2} - \frac{rmX^*Y^*Z^*}{K(a+Y)^2}\bigg) \nonumber \\
	& + & \bigg(\frac{s^\alpha}{\alpha\Gamma (\alpha)}\bigg)^2 X^*Y^* \bigg(\frac{r\lambda}{K} + \lambda^2\bigg), \nonumber\\
	A_3 & = & 1 +  \frac{s^\alpha}{\alpha\Gamma (\alpha)}\bigg(\frac{mY^*Z^*}{(a+Y^*)^2} - \frac{rX^*}{K}\bigg) + \bigg(\frac{s^\alpha}{\alpha\Gamma (\alpha)}\bigg)^2 \frac{amdZ^*}{(a+Y^*)^2} \bigg(1 - \frac{s^\alpha}{\alpha\Gamma (\alpha)}\frac{rX^*}{K}\bigg) \nonumber \\
	& + & \bigg(\frac{s^\alpha}{\alpha\Gamma (\alpha)}\bigg)^2 \bigg(\frac{rmX^*Y^*Z^*}{K(a+Y^*)^2} + X^*Y^* \bigg(\frac{r\lambda}{K} + \lambda^2\bigg) \bigg) \nonumber.
	\end{eqnarray} 
	
	\noindent We already know that the interior equilibrium $E^{*}$ exists if $(i)$ $R_{0} >1$ and $(ii)$ $\theta > \theta_{1}$, where $\theta_1 = d + \frac{\lambda ad(r + \lambda K)}{r(\lambda K - \mu)}$. Using these conditions, we prove the conditions of Lemma $2$:
	\begin{itemize}
		\item[(i)] $p(1) = 1 + A_1 +A_2 + A_3 \\
		= \bigg(\frac{s^\alpha}{\alpha\Gamma (\alpha)}\bigg)^3 \bigg(\frac{adrmX^*Z^*}{K(a + Y^*)^2}\bigg)$	$> 0$  for any $s$.
		
		\item[(ii)] $(-1)^3p(-1) = 1 - A_1 +A_2 - A_3 \\
		= 4\bigg(2 - \frac{s^\alpha}{\alpha\Gamma (\alpha)} \frac{rX^*}{K}\bigg) + \frac{s^\alpha}{\alpha\Gamma (\alpha)}\frac{mZ^*}{(a+I^*)^2}\bigg(\frac{s^\alpha}{\alpha\Gamma (\alpha)}ad+2Y^*\bigg)\bigg(2 - \frac{s^\alpha}{\alpha\Gamma (\alpha)}\frac{rX^*}{K}\bigg)\\
		+ \bigg(\frac{s^\alpha}{\alpha\Gamma (\alpha)}\bigg)^2 2X^*Y^*\bigg(\frac{r\lambda}{K} + \lambda^2\bigg)$	$> 0$, 
		if $s < \sqrt[\alpha]{\frac{2K\alpha\Gamma(\alpha)}{rX^*}}$.

		\item[(iii)] $(1 - A_3^2) - \mid A_2 - A_3 A_1 \mid = \bigg[ 1 -  \bigg\{\bigg(1 + \frac{s^\alpha}{\alpha\Gamma (\alpha)}\frac{mZ^*}{(a+Y^*)^2} \bigg(\frac{s^\alpha}{\alpha\Gamma (\alpha)}ad + Y^*\bigg)\bigg)\bigg(1 - \frac{s^\alpha}{\alpha\Gamma (\alpha)}\frac{rX^*}{K}\bigg) + \bigg(\frac{s^\alpha}{\alpha\Gamma (\alpha)}\bigg)^2 X^*Y^*\bigg(\frac{r\lambda}{K} + \lambda^2\bigg)\bigg\}^2 \bigg] - \\
		\left| \bigg\{3 +  \frac{2s^\alpha}{\alpha\Gamma (\alpha)}\bigg(\frac{mY^*Z^*}{(a+Y^*)^2} - \frac{rX^*}{K}\bigg) + \bigg(\frac{s^\alpha}{\alpha\Gamma (\alpha)}\bigg)^2 \bigg(\frac{amdZ^*}{(a+Y^*)^2} - \frac{rmX^*Y^*Z^*}{K(a+Y^*)^2}\bigg)\right.$\\
		$\left. + \bigg(\frac{s^\alpha}{\alpha\Gamma (\alpha)}\bigg)^2 X^*Y^* \bigg(\frac{r\lambda}{K} + \lambda^2\bigg)\bigg\} - \bigg(\frac{s^\alpha}{\alpha\Gamma (\alpha)}\frac{rX^*}{K} - \frac{s^\alpha}{\alpha\Gamma (\alpha)} \frac{mY^*Z^*}{(a+Y^*)^2} - 3\bigg)\right.$\\
		$\left.\bigg\{\bigg(1 + \frac{s^\alpha}{\alpha\Gamma (\alpha)}\frac{mZ^*}{(a+Y^*)^2} \bigg(\frac{s^\alpha}{\alpha\Gamma (\alpha)}ad + Y^*\bigg)\bigg)\bigg(1 - \frac{s^\alpha}{\alpha\Gamma (\alpha)}\frac{rX^*}{K}\bigg) + \right.$\\ $\left.\bigg(\frac{s^\alpha}{\alpha\Gamma (\alpha)}\bigg)^2 X^*Y^*\bigg(\frac{r\lambda}{K} + \lambda^2\bigg)\bigg\}\right|$. 
	\end{itemize}
	Due to the complexity of mathematical expression, it is difficult to find a restriction on the step-size $s$ in terms of other system parameters so that $(1 - A_3^2) - \mid A_2 - A_3 A_1 \mid > 0$. However, in the numerical section, we show that there exists a value  $s_9$ (say) of $s$ such that 
	\begin{equation}\label{Jury}
	(1 - A_3^2) - \mid A_2 - A_3 A_1 \mid >0, ~\forall s < s_9.
	\end{equation}
	Thus, following Jury criteria, the fixed point $E^*$ is locally asymptotically stable for $0<\alpha\leq 1$ if $R_0 > 1$, $\theta > \theta_{1}$ and $0 < s < min\{s_{8}, s_{9}\}$, where  $s_8 = \sqrt[\alpha]{\frac{2K\alpha\Gamma(\alpha)}{rX^*}}$,  $\theta_1=d + \frac{\lambda ad(r + \lambda K)}{r(\lambda K - \mu)}$ and $s_{9}$ is such that \eqref{Jury} holds. Hence the theorem.

\section{Numerical Simulations}

In this section, an extensive numerical simulations are performed for the fractional-order discrete system (\ref{Discrete model}) to validate our theoretical results for different fractional values of $\alpha$ $(0 < \alpha < 1)$. Stability of fixed points depends on the step-size $s$ and fractional-order $\alpha$ (see Theorem 1). To demonstrate the effect of fractional order $\alpha$ and step-size $s$ on the system dynamics, we cite different examples.\\

\noindent{\bf Table 3.1.} Restrictions on the step-size, following Theorem 1(b,c), for the stability of fixed points $E_1$ and $E_2$ for different fractional-orders $\alpha$. For $E_1$, we consider  $K = 40$, $\lambda = 0.005$, $\theta = 0.189$ and for $E_2$ we consider $K = 200$, $\lambda = 0.015$, $\theta = 0.08$. Following parameters are common for both the equilibrium points: $r = 2.0$, $m = 0.52$, $\mu = 0.28$, $a = 15.0$, $d = 0.09$. \\

{\small
	\begin{tabular} 
		{|l l l|} \hline
		& $E_1$ & $E_2$ \\ \hline
		Fractional order $\alpha$ & Step-size $s<min(s_2, s_3, s_4)$ & Step-size $s_5< s <min(s_6,s_7)$  \\  \hline
		
		$\alpha = 0.3$                 &  $s_2 = 21513$         & $s_5 = 0.0248$ \\
		&  $s_3 = 0.6973$           & $s_6 = 1835600$  \\
		&  $s_4 = 31857$           & $s_7 = 2.1578$  \\
		\hline
		$\alpha = 0.4$                 & $s_2 = 1726.2$          & $s_5 = 0.0608$ \\
		&  $s_3 = 0.7415$           & $s_6 = 48464$   \\
		&  $s_4 = 2317.3$           & $s_7 = 1.7302$  \\
		\hline
		$\alpha = 0.6$                 & $s_2 = 145.5959$           & $s_5 = 0.1564$ \\
		&  $s_3 = 0.8289$         & $s_6 = 1344.9$   \\
		&  $s_4 = 177.1754$           & $s_7 = 1.4582$  \\
		\hline
		$\alpha = 0.8$                 & $s_2 = 44.1464$            & $s_5 = 0.2619$ \\
		&   $s_3 = 0.9150$          & $s_6 = 233.9135$   \\
		&  $s_4 = 51.1488$           & $s_7 = 1.3976$  \\
		\hline
		$\alpha = 0.95$                & $s_2 = 25.6082$              & $s_5 = 0.3414$ \\
		&  $s_3 = 0.9788$            & $s_6 = 104.2794$   \\
		&  $s_4 = 28.9884$           & $s_7 = 1.3984$  \\
		\hline
\end{tabular}}
\label{Table 3.1}
\vspace{1cm}

\begin{figure}[H]
	%\hspace{2in}
	\includegraphics[width=9in, height=2.5in]{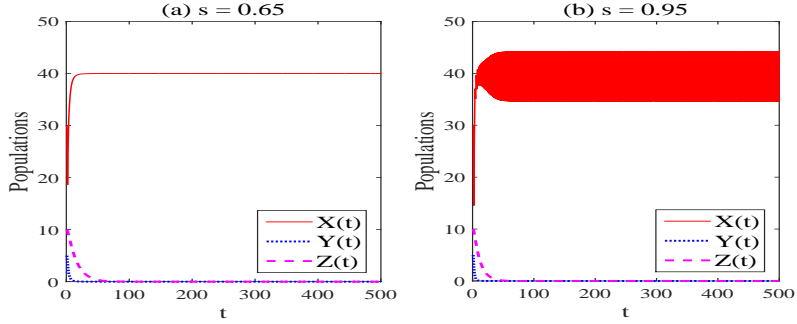}
	\vspace{-2.5cm}
	\caption{Stable (Fig. a) and unstable (Fig. b) behavior of the fixed point $E_1$ for the step-size $s = 0.65$ and $s = 0.95$ corresponding to the fraction order $\alpha=0.8$. Parameter values are $r = 2.0$, $K = 40.0$, $\lambda = 0.005$, $m = 0.52$, $\mu = 0.28$, $a = 15.0$, $\theta = 0.189$, $d = 0.09$. \label{limit_cycle_2}}
\end{figure}

\noindent \textbf{Example 1:}
We consider the following parameter values from \cite{Bairagi08}: $r = 2.0$, $K = 40.0$, $\lambda = 0.005$, $m = 0.52$, $\mu = 0.28$, $a = 15.0$, $\theta = 0.189$, $d = 0.09$ and initial point $S(0) = 30, I(0) = 5, P(0) = 10$. With these parameter values, we calculate different limits $s_i$ on the step-size for the stability of fixed points $E_1$ and $E_2$ for different values of fractional-order $\alpha$. For example, we compute $s_2=44.1464$, $s_3=0.9150$ and $s_4=51.1488$ for $\alpha=0.8$. The equilibrium $E_1$ will then be locally asymptotically stable, following Theorem 1(b), if $s< min(44.1464,0.9150, 51.1488)= 0.9150$. It is unstable if $s>0.9150$. Similarly, $E_2$ will be locally asymptotically stable, following Theorem 1(c), if $0.2619<s< min(233.9135, 1.3976)= 1.3976$. Similar step-size restrictions for the stability of the above two equilibrium points for different fractional order $\alpha$ are tabulated in Table 3.1. Considering $s=0.65 (<0.9150)$ and $s=0.95 (>0.9150)$ corresponding to $\alpha=0.8$, the behavior of the solution trajectories of the system (\ref{Eco-epidemiological fractional order model}) are plotted in Fig. \ref{limit_cycle_2}. It shows that the predator-free equilibrium $E_1$ is stable for the step-size $s=0.65$ and unstable for $s=0.95$. 

%\newpage 

\begin{figure}[H]
	\includegraphics[width=9in, height=2.5in]{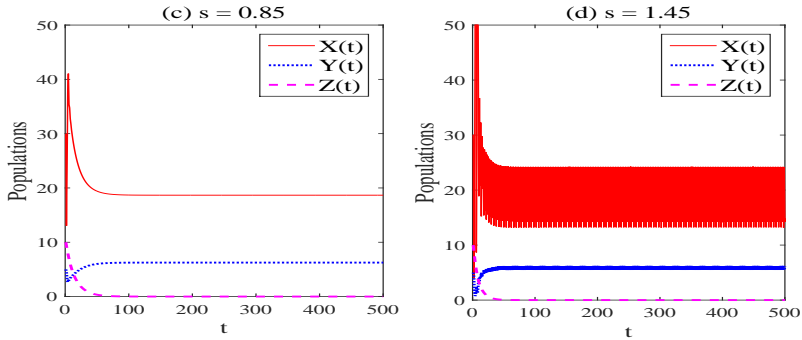}
	\vspace{-2cm}
	\caption{Stable (Fig. a) and unstable (Fig. b) behavior of the fixed point $E_2$ for the step-size $s = 0.85$ and $s = 1.45$ corresponding to the fraction order $\alpha=0.8$. Parameter values are $\lambda = 0.005$, $K = 200$, $\theta = 0.08$. Other parameters and initial values are as in  Fig. \ref{limit_cycle_2}.	\label{limit_cycle_3}}
\end{figure}
\noindent Similarly, we observe (Figure \ref{limit_cycle_3}) that the predator-free fixed point $E_2$ of the system (\ref{Eco-epidemiological fractional order model}) is stable for $s=0.85$ and unstable for $s=1.45$ corresponding to order $\alpha=0.8$ (see Theorem 1(c) and Table 3.1.).
\begin{figure}[H]
	%	\hspace{1cm}
	\includegraphics[width=6in, height=3.0in]{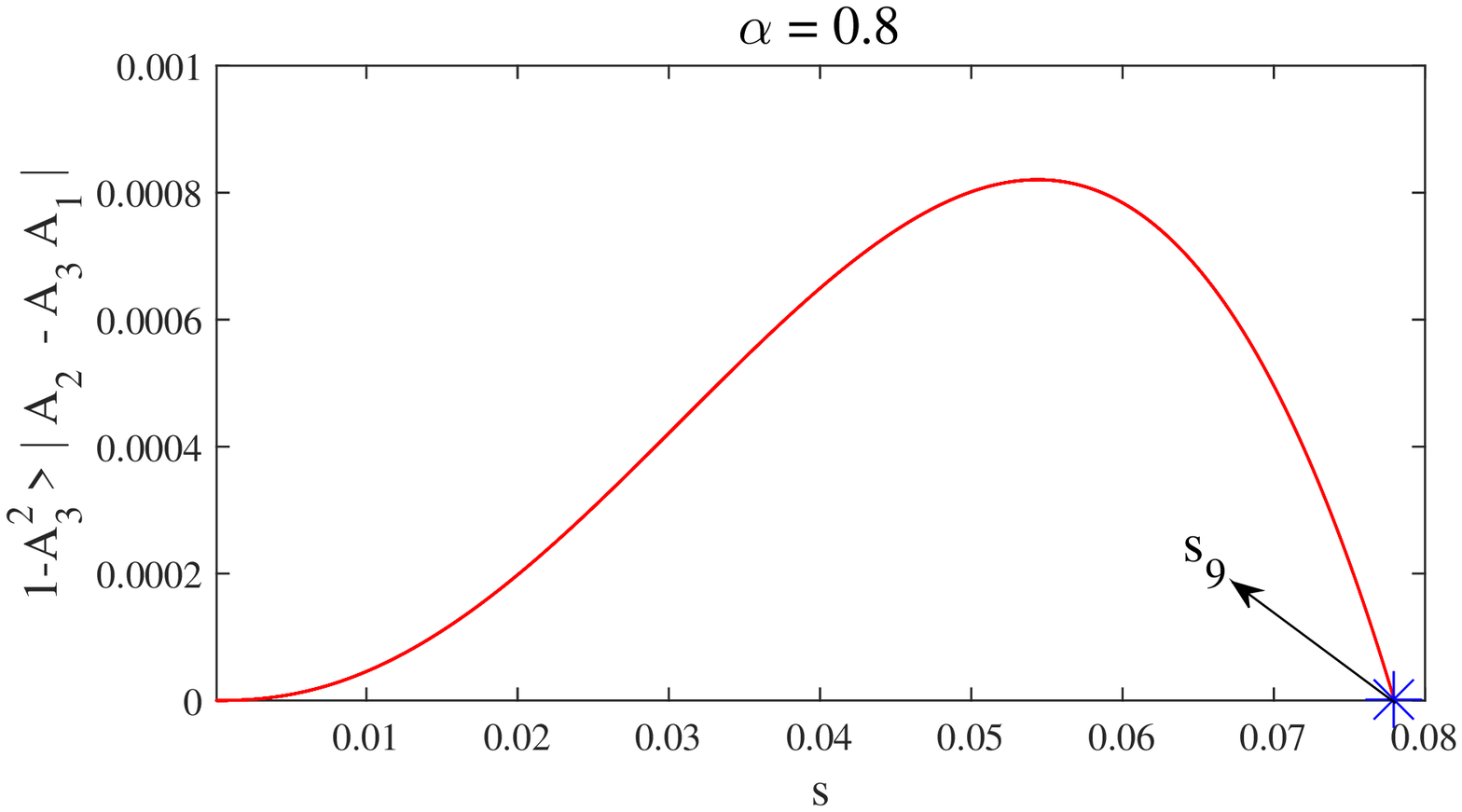}
	\vspace{-3.5cm}
	\caption{Existence of the bound $s_9$ such that the condition $1-A_3^2 >$ $ \mid A_2 - A_3 A_1 \mid$ is hold. The bound is obtained as $s_9 = 0.07798$ for the parameter values $r = 15.0$, $K = 40.0$, $\lambda = 0.006$, $m = 14.5$, $\mu = 0.0019$, $a = 16.0$, $\theta = 11.1$, $d = 6.0$ and $\alpha=0.8$.}
	\label{limit_cycle_4}
\end{figure}
\noindent{\bf Table 3.2.} Restriction on the step-size, following Theorem 1(d), for the stability of fixed point $E^*$ for different fractional-order $\alpha$. Parameter values are as in Fig. \ref{limit_cycle_4}. \\

\hspace{1in}
{\small
	\begin{tabular}
		{|l l|} \hline
		& $E^*$  \\ \hline
		Fractional order $\alpha$ & Step-size $s< min(s_8, s_9)$   \\  \hline
		
		$\alpha = 0.3$                 &  $s_8 = 0.0074$          \\
		& $s_9 = 0.00098$  \\
		\hline
		$\alpha = 0.4$                  &  $s_8 = 0.0245$          \\
		& $s_9 = 0.00538$  \\
		\hline
		$\alpha = 0.6$                 &  $s_8 = 0.0853$          \\
		& $s_9 = 0.03108$  \\
		\hline
		$\alpha = 0.8$                  &  $s_8 = 0.1662$          \\
		& $s_9 = 0.07798$  \\
		\hline
		$\alpha = 0.95$                 &  $s_8 = 0.2327$          \\
		& $s_9 = 0.1073$  \\
		\hline
\end{tabular}}
\vspace{0.5cm}

\noindent \textbf{Example 2:}
We consider another parameter set as in \cite{ChattoBairagi01}: $r = 15.0$, $K = 40.0$, $\lambda = 0.006$, $m = 14.5$, $\mu = 0.0019$, $a = 16.0$, $\theta = 11.1$, $d = 6.0$ and initial point $S(0) = 30, I(0) = 5, P(0) = 10$. First we numerically determine the bound $s_9=0.07798$ (Fig. \ref{limit_cycle_4}) for the stability of the interior fixed point $E^*$
corresponding to a fractional order $\alpha=0.8$ so that the condition $1-A_3^2 > \mid A_2 - A_3 A_1 \mid$ is satisfied (see Theorem 1(d)).

\noindent As before, we construct Table 3.2 to show the bounds on step-size for the stability of the interior fixed point $E^*$ for different values of fractional-order $\alpha$. One can also compute that $R_0 = 126.3158 > 1$ and $\theta =11.1> \theta_{1} = 8.4579$ for $\alpha=0.8$. Following Theorem 1(d), the fixed point $E^* (20.8573, 18.8235,$ $0.2962)$ is then stable (Fig. 4a) for $s = 0.05< min(s_8, s_9)=0.07798$ (see Table 3.2) and unstable for $s = 0.08 (>0.07798)$ (Fig. 4b) for the fractional-order $\alpha = 0.8$.

\begin{figure}[H]
	%	\hspace{2in}
	\includegraphics[width=9in, height=3.0in]{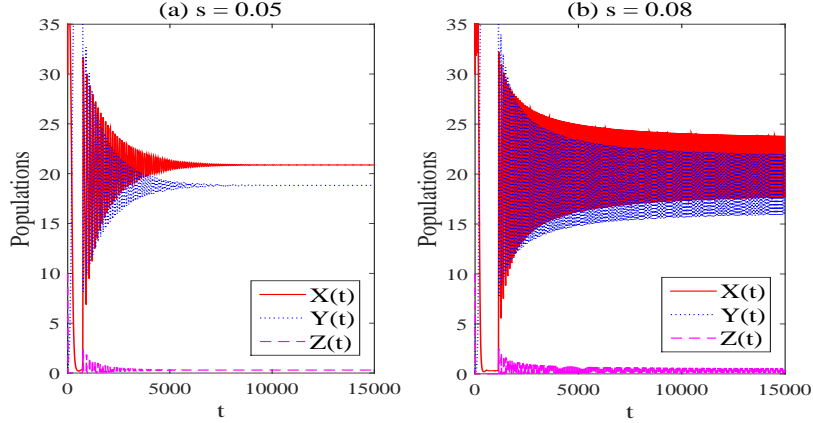}
	\vspace{-2cm}
	\caption{Stable and unstable behavior of the fixed point $E^*$, respectively, for the step-sizes $s = 0.05$ (Fig. a) and $s = 0.08$  (Fig. b) when $\alpha = 0.8$. Here all parameters are as in Figure \ref{limit_cycle_4}.}
	\label{limit_cycle_5}
\end{figure}
\vspace{-0.75cm}
\begin{figure}[H]
	\hspace{-0.5cm}
	\includegraphics[width=9in, height=3.0in]{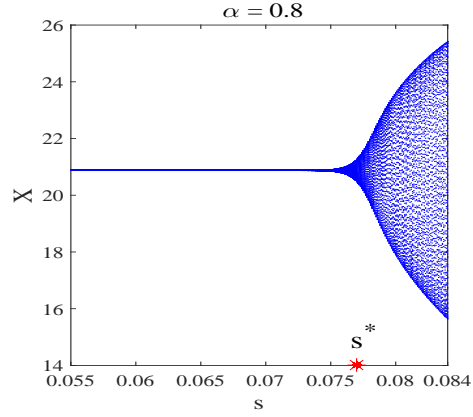}
	\vspace{-2.5cm}
	\caption{Bifurcation diagram of $X$ population of system (\ref{Discrete model}) in the range $[0.055, 0.085]$ with step-size $s$ as the bifurcation parameter. System remains stable for  $s < s^* $ and becomes unstable for $s > s^* $, where $s^* =  0.077$. Here $\alpha = 0.8$ and othere parameters are as in Figure \ref{limit_cycle_4}.}
	\label{limit_cycle_6}
\end{figure}
\noindent Bifurcation diagram of $X$ population (Fig. \ref{limit_cycle_6}) shows that population is stable for $s<s^*$ and unstable for $s>s^*$, where $s^*=0.077$. These results show that stability of the discrete fractional order system (\ref{Eco-epidemiological fractional order model}) strongly depends on the step-size and the fractional order.\\

\begin{figure}[H]
	\includegraphics[width=9in, height=2.5in]{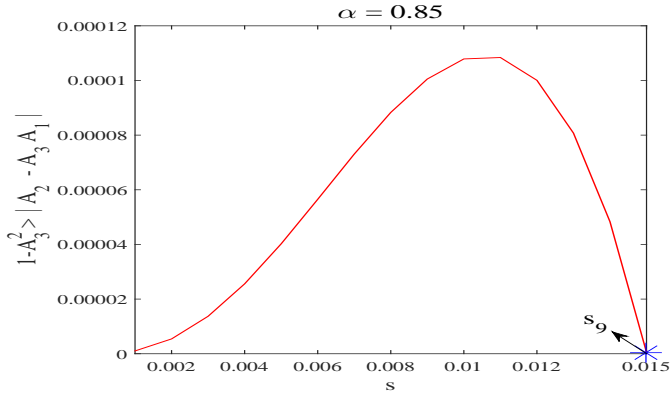}
	\vspace{-1.5cm}	
	\caption{The bound $s_9 = 0.0150$ is determined such that the condition $1-A_3^2 >$ $ \mid A_2 - A_3 A_1 \mid$ holds. The parameter values are $r = 22.0$, $K = 300.0$, $\lambda = 0.06$, $m = 15.5$, $\mu = 2.3$, $a = 15.0$, $\theta = 10$, $d = 8.3$ and $\alpha=0.85$.} 
	\label{limit_cycle_7}
\end{figure}

\begin{figure}[H]
	\includegraphics[width=8in, height=2.5in]{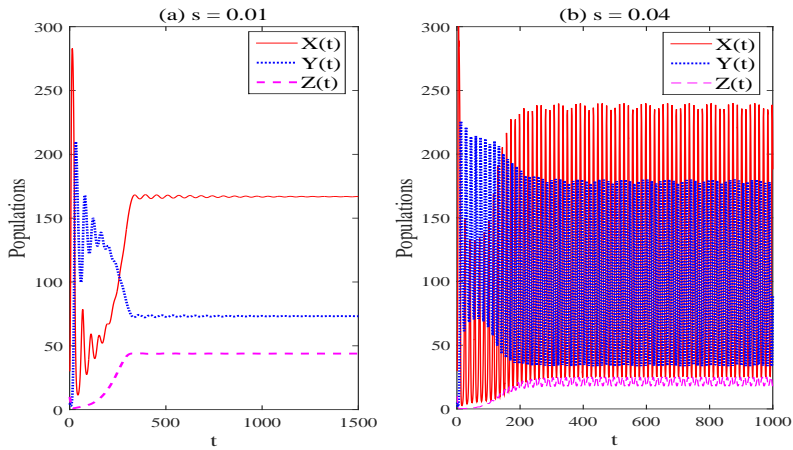}
	\caption{Stable and unstable behavior of the fixed point $E^* = (166.8449, $ $73.2353,$ $43.8939)$ for the step-sizes $s = 0.01$ (Fig. a) and $s = 0.04$ (Fig. b) when $\alpha = 0.85$. Here all parameters are as in Fig. \ref{limit_cycle_7}.}
	\label{limit_cycle_8}
\end{figure}

\noindent \textbf{Example 3:}
In \cite{ChattoBairagi08}, the authors showed that the continuous system may show chaotic dynamics for some parameter values. Here we consider the same parameter set and initial value for which the system \eqref{Eco-epidemiological integer order model} shows chaotic behavior: $r = 22.0$, $K = 300.0$, $\lambda = 0.06$, $m = 15.5$, $\mu = 2.3$, $a = 15.0$, $\theta = 10$, $d = 8.3$ and $(S(0), I(0), P(0))=(30, 5, 10)$. As before, we determine the bound $s_9=0.0150$ (see Fig. \ref{limit_cycle_7}) for the stability of fixed point $E^*$ when $\alpha=0.85$ and tabulate (see Table 3.3) the bounds on step-size, following Theorem 1(d), for different values of fractional-order $\alpha$. Table 3.3 shows that the critical value of the step-size, where the switching of stability occurs, decreases as the order of the fractional derivative decreases.
In Figure \ref{limit_cycle_8}, we present the stable (Fig. a) and unstable (Fig. b) behavior of the system around the interior fixed point $E^* = (166.8449, $ $73.2353,$ $43.8939)$ for the step-sizes $s = 0.01$ and $s = 0.04$, respectively, when $\alpha = 0.85$. 

\noindent{\bf Table 3.3.} Restriction on the step-size for the stability of fixed point $E^*$ for different fractional-orders $\alpha$ with the parameters of Fig. \ref{limit_cycle_7} .\\

\hspace{1in}
{\small
	\begin{tabular}
		{|l l|} \hline
		& $E^*$  \\ \hline
		Fractional order $\alpha$ & Step-size $s< min(s_8, s_9)$   \\  \hline
		$\alpha = 0.95$                 &  $s_8 = 0.1455$,  $s_9 = 0.0242$         \\
		%	& $s_9 = 0.0242$  \\
		\hline
		$\alpha = 0.85$                  &  $s_8 = 0.1112$, $s_9 = 0.0150$           \\
		%	& $s_9 = 0.0150$  \\
		\hline
		$\alpha = 0.6$                 &  $s_8 = 0.0405$, $s_9 = 0.0023$          \\
		%& $s_9 = 0.0023$  \\
		\hline
		$\alpha = 0.55$                  &  $s_8 = 0.0300$ , $s_9 = 0.0013$          \\
		%	& $s_9 = 0.0013$  \\
		\hline
		
		$\alpha = 0.45$                 &  $s_8 = 0.0136$, $s_9 = 0.0003$           \\
		%	& $s_9 = 0.0003$  \\
		\hline
\end{tabular}}

%\newpage

\begin{figure}[H]
	\includegraphics[width=7.5in, height=2.5in]{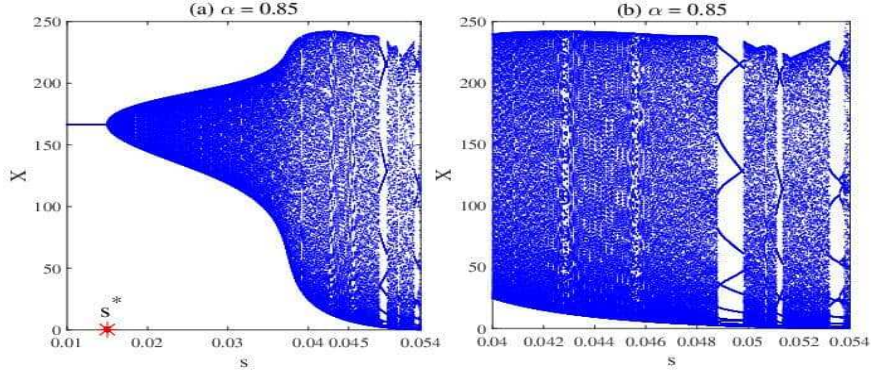}
	\vspace{-1.25cm}
	\caption{Bifurcation diagrams of $X$ population of system (\ref{Discrete model}) in the range $[0.01, 0.054]$ with respect to the step-size $s$. System is stable for $s <s^*$ and unstable for  $s > s^* $. Exchange of stability occurs at $s^*=  0.015$ (Fig. 8a). A magnified view of the bifurcation diagram for $0.04\leq s \leq 0.054$ is presented in Fig. 8b. All parameters are as in Fig. 6 with $\alpha=0.85$.}
	\label{limit_cycle_9}
\end{figure}

\noindent Bifurcation diagram of system populations (Fig. 8a) shows that population is stable for $s<s^*$ and losts its stability for $s>s^*$ corresponding to fractional order $\alpha = 0.85$, where $s^*=0.015$. As the step-size is further increased, the system shows complex dynamics like chaos. An enlarged view (Fig. 8b)  of the bifurcation diagram for small range of higher step-size is presented for visualisation of the chaotic dynamics. 
Existence of such chaotic dynamics (see Fig. 9) may be verified by computing the largest Lyapunov exponent \cite{AET18}.
\begin{figure}[H]
	%\hspace{-0.5cm}
	\includegraphics[width=8in, height=3in]{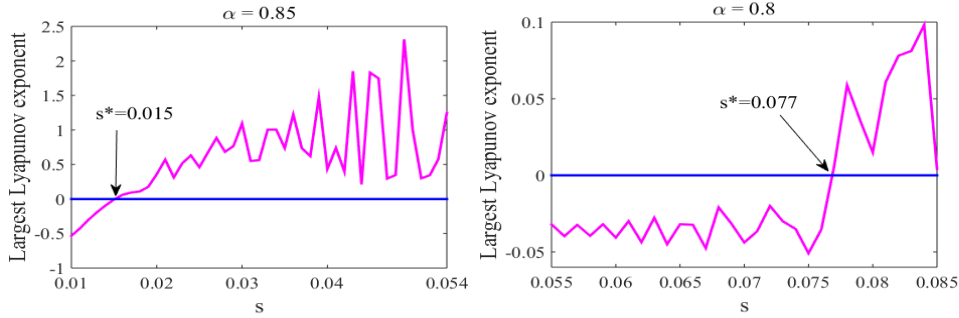}
	\vspace{-3.0cm}
	\caption{Largest Lyapunov exponent (LLE) drawn for two values of $\alpha$. The left figure shows that LLE becomes positive for $s>s^*$, where $s^*=0.015$, $\alpha=0.85$ and other parameters are as in Fig. 8. The right figure shows that LLE becomes positive for $s>s^*$, where $s^*=0.077$, $\alpha=0.8$ and other parameters are as in Fig. 5.} 
	%	\label{lyapunov}
\end{figure}
\section{Summary}

The application of fractional calculus in the study of biological systems is important because it considers past history of the state variables. On the other hand, discretization of a continuous system is inevitable for numerical computations and for population models, particularly where species generation does not overlap. It is earlier shown that the stability of the discrete fractional-order system depends on its fractional-order \cite{ Matouk15, Alzabut18, SJ18, SET18,Salman17,AET18}. It was shown for a simpler two-dimensional epidemic model that the stability of a fractional-order discrete system also depends on the step-size \cite{Raheem14}. It is worth mentioning that the determination of explicit bounds would be difficult if the nonlinearity and dimension of the system are high. The novelty of this work is that we determined the analytical bounds on the step-size for a highly nonlinear three-dimensional system so that the stability of an equilibrium point is maintained. We considered a highly nonlinear three-dimensional system of fractional-order differential equations that represent the predator-prey interaction in the presence of infection. We revealed the dynamics of the corresponding discrete fractional system. Discretization of the fractional-order system was done with piecewise constant arguments and the corresponding dynamics were explored around different feasible fixed points. We used the Jury criterion for local stability of the discrete fractional-order systems. It is observed that stability of the system depends on both the step-size and fractional order. Different examples were cited to illustrate the stability of predator-free, infection-free, and coexistence equilibrium points. These simulation results perfectly agree with the analytical results. More specifically, the critical value of the step-size, where the switching of stability occurs, decreases as the order of the fractional derivative decreases. Bifurcation diagrams with respect to step size show that the system is stable and all populations coexist in the stable state if the step-size does not exceed some threshold value, which is also dependent on the fractional order, and unstable if it exceeds. Our simulation result also shows that the system may show chaotic dynamics, as shown in the continuous-time case, if the step-size is large. This study thus reveals that the dynamics of a discrete fractional order system strongly depends on the step-size and the fractional order.
\\

%\noindent {\bf Acknowledgement:} We are thankful to the reviewers for their suggestions.  

%%%%%%%% Bibliography %%%%%%%%%%%%%%%%%%%%%%%%%%%%%%%%%%%%%%%%%%%%%%%%%


\begin{thebibliography}{20}
	
	
	\bibitem{SET93}
	S. Samko, A. Kilbas, O. Marichev,
	\emph{Fractional Integrals and Derivatives: Theory and Applications}.
	Gordon and Breach Science Publishers, Yverdon, (1993).
	
	\bibitem{Rihan15}
	F.A. Rihan, S. Lakshmanan, A. H. Hashish, R. Rakkiyappan, E. Ahmed,
	Fractional-order delayed predator–prey systems with Holling type-II functional response.
	\emph{Nonlinear Dynamics}.
	\textbf{80}, 1-2 (2015), 777--789.
	
	
	\bibitem{MET11} 
	J. T. Machado, V. Kiryakova, F. Mainardi,
	Recent history of fractional calculus.
	\emph{Commun. Nonlinear Sci. Numer. Simulat.} \textbf{ 16}, (2011), 1140--1153.
	
	\bibitem{MO15} 
	I. Matychyn, V. Onyshchenko, Time-optimal control of fractional-order linear systems. 
	\emph{Fract. Calc. Appl. Anal.} \textbf{18},  3 (2015), 687--696.
	
	\bibitem{GET10} 
	R. E. Gutierrez, J. M. Rosario, J. T. Machado,
	Fractional Order Calculus: Basic Concepts and Engineering Applications.
	\emph{Math. Prob. Eng.} (2010); DOI: 10.1155/2010/375858.
	
	\bibitem{AET12} 
	S. Abbas, M. Benchohra, G. M. N'Guerekata,
	\emph{Topics in fractional differential equations}.
	Springer Science \& Business Media, (2012).
	
	\bibitem{D11} 
	S. Das,
	\emph{Introduction to fractional calculus for scientists and engineers}.
	Springer, (2011).
	
	\bibitem{M06} 
	M. L. Richard,
	\emph{Fractional calculus in bioengineering}.
	Redding: Begell House, (2006).
	
	\bibitem{LZ18} 
	Y. Li and Q. Zhang, Blow-up and global existence of solutions for a time fractional diffusion equation.
	\emph{Fract. Calc. Appl. Anal.} \textbf{21},  6 (2018), 1619--1640.
	
	\bibitem{M13} 
	J. Munkhammar, Chaos in a fractional order logistic map,
	\emph{Fract. Calc. Appl. Anal.} \textbf{16},  3 (2013), 511--519.
	
	\bibitem{TET15} 
	M. C. Tripathy, D. Mondal, K. Biswas, S. Sen,
	Experimental studies on realization of fractional inductors and fractional order bandpass filters.
	\emph{Int. J. Circuit Theory and Applications.} \textbf{43}, 9 (2015), 1183--1196.
	
	\bibitem{TB84} 
	P. J. Torvik, R. L. Bagley,
	On the appearance of the fractional derivative in the behaviour of real materials.
	\emph{J. Appl. Mechanics.} \textbf{51},  2 (1984), 294--298.
	
	\bibitem{SET11} 
	J. A. Sabatier, O. P. Agrawal, J. T. Machado,
	\emph{ Advances in fractional calculus}.
	Dordrecht: Springer, (2007).
	
	\bibitem{BET17}
	A. Boukhouima, K. Hattaf, N. Yousfi,
	Dynamics of a fractional order HIV infection model with specific functional response and cure rate.
	\emph{Int. J. Diff. Equations.} \textbf{2017}, (2017).
	
	
	\bibitem{Ahmed07} 
	E. Ahmed, A. M. A. El-Sayed, H. A. A. El-Saka,
	Equilibrium points, stability and numerical solutions of fractional-order predator-prey and rabies models.
	\emph{J. Math. Anal. Appl.} \textbf{325}, (2007), 542--553.
	
	\bibitem{RET13} 
	S. Ranaa, S. Bhattacharyaa, J. Pal, G. M. N'Guerekata, J. Chattopadhyay,
	Paradox of enrichment: A fractional differential approach with memory.
	\emph{Physica A.} \textbf{392}, (2013), 3610--3621.
	
	\bibitem{CY14} 
	Z. Cui, Z. Yang,
	Homotopy perturbation method applied to the solution of fractional lotka-volterra equations with variable coefficients.
	\emph{J. Mod. Meth. Numer. Math.} \textbf{5}, (2014), 1--9.
	
	\bibitem{MET19} 
	S. Mondal, N. Bairagi, G. M. N'Guerekata,
	Global stability of a Leslie-Gower-type fractional order tritrophic food chain model.
	\emph{Fractional Differential Calculus.} \textbf{9}, 1 (2019), 149--161. 
	
	\bibitem{MET17} 
	S. Mondal, N. Bairagi, A. Lahiri,
	A fractional calculus approach to Rosenzweig-MacArthur predator-prey model and its solution.
	\emph{J. Mod. Meth. Numer. Math.} \textbf{8}, 1-2 (2017), 66--76.
	
	\bibitem{HET15} 
	J. Huo, H. Zhao, L. Zhu,
	The effect of vaccines on backward bifurcation in a fractional order HIV model.
	\emph{Nonlinear Anal. RWA.} \textbf{26}, (2015), 289--305.
	
	\bibitem{LET16} 
	H. L. Li, L. Zhang, C. Hu, Y. L. Jiang, Z. Teng,
	Dynamical analysis of a fractional-order predator-prey model incorporating a prey refuge.
	\emph{J. Appl. Math. Comput.} (2016); DOI: 10.1007/s12190-016-1017-8.
	
	\bibitem{V15} 
	C. Vargas-De-Leon,
	Volterra-type Lyapunov functions for fractional-order epidemic systems.
	\emph{Commun. Nonlinear Sci. Numer. Simul.} \textbf{24}, (2015), 75--85.
	
	
	\bibitem{MET17a} 
	S. Mondal, A. Lahiri, N. Bairagi,
	Analysis of a fractional order eco-epidemiological model with prey infection and type 2 functional response.
	\emph{Math. Meth. Appl. Sci.} \textbf{40},  18 (2017), 6776--6789.
	
	\bibitem{Hattaf15}
	K. Hattaf, A. A. Lashari, B. El Boukari,N. Yousfi,
	Effect of discretization on dynamical behavior in an epidemiological model.
	\emph{Differential Equations and Dynamical Systems.} \textbf{23}, 4 (2015), 403--413.
	
	\bibitem{Sekiguchi10}
	M. Sekiguchi, E. Ishiwata,
	Global dynamics of a discretized SIRS epidemic model with time delay.
	\emph{Journal of Mathematical Analysis and Applications.} \textbf{371}, 1 (2010), 195--202.
	
	\bibitem{Chen13}
	Q. Chen, Z. Teng, L. Wang, H. Jiang,
	The existence of codimension-two bifurcation in a discrete SIS epidemic model with standard incidence.
	\emph{Nonlinear Dyn.} \textbf{71},  (2013), 55--73.
	
	\bibitem{Hu12}
	Z. Hu, Z. Teng, L. Wang, H. Jiang,
	Stability analysis in a class of discrete SIRS epidemic models.
	\emph{Nonlinear Anal.: Real World Appl. } \textbf{13},  (2012), 2017--2033.
	
	\bibitem{Biswas17}
	M. Biswas, N. Bairagi,
	Discretization of an eco-epidemiological model and its dynamic consistency.
	\emph{Journal of Difference Equations and Applications.} \textbf{23}, 5 (2017), 860--877.
	
	\bibitem{GET20}
	S. Ghorai, P. Chakraborty, S. Poria and N. Bairagi,
	Dispersal-induced pattern forming instabilities in host-parasitoid metapopulations.
	\emph{Nonlinearity.} (2020), DOI.org/10.1007/s11071-020-05505-w.
	
	
	\bibitem{AE07} 
	F. M. Atici, P. W. Eloe,
	A transform method in discrete fractional calculus.
	\emph{Int. J. Difference Equations.} \textbf{2}, 2 (2007), 165--176.
	
	\bibitem{CET1107} 
	C. Fulai, L. Xiannan, Z. Yong,
	Existence results for nonlinear fractional difference equation.
	\emph{Advances in Difference Equations.} \textbf{2011}, (2011); DOI: 10.1155/2011/713201.
	
	\bibitem{WB14} 
	G. C. Wu, D. Baleanu,
	Chaos synchronization of the discrete fractional logistic map.
	\emph{Signal Processing.} \textbf{102}, (2014), 96--99.
	
	\bibitem{Matouk15} 
	A. A. Elsadany, A. E. Matouk,
	Dynamical behaviors of fractional-order Lotka Volterra predator-prey model and its discretization.
	\emph{J. Appl. Math. Comput.} \textbf{49}, (2015), 269--283.
	
	\bibitem{Alzabut18} 
	J. Alzabut	, T. Abdeljawad, D. Baleanu,
	Nonlinear delay fractional difference equations with applications on discrete fractional Lotka-Volterra competition model.
	\emph{J. Comput. Anal. Appl.} \textbf{25},  5 (2018), 889--898.
	
	\bibitem{SJ18} 
	A.G. M. Selvam, R. Janagaraj,
	Numerical Analysis of a Fractional Order Discrete Prey – Predator System with Functional Response.
	\emph{In. J. Engineering \& Technology.} \textbf{7},  4.10 (2018), 681--684.
	
	\bibitem{SET18} Y Shi, Q. Ma and X. Ding, 2018. Dynamical behaviors in a discrete fractional-order predator-prey system. \emph{Filomat.} \textbf{32}, 17 (2018), 5857--5874.
	
	\bibitem{Salman17} 
	S. M. Salman,
	On a Discretized Fractional-Order SIR Model for Influenza.
	\emph{Progress in Fractional Differentiation and Application.} \textbf{3},  2 (2017), 163--173. 
	
	\bibitem{AET18} A. M. Abdelaziz, A. I. Ismail, F. A. Abdullah and M. H. Mohd, 2018. Bifurcations and chaos in a discrete SI epidemic model with fractional order. \emph{Advances in Difference Equations.}  1 (2018), 1--19.
	
	
	\bibitem{Khan19} 
	A. Khan, J. F. Gomez-Aguilar, T. S. Khan, H. Khan,
	Stability analysis and numerical solutions of fractional order HIV/AIDS model.
	\emph{Chaos, Solitons and Fractals.} \textbf{122}, (2019), 119--128. \\
	
	\bibitem{Khan20}
	A. Khan, T. Abdeljawad, J. F. Gomez-Aguilar, H. Khan,
	Dynamical study of fractional order mutualism parasitism food web module. 
	\emph{Chaos, Solitons and Fractals.} \textbf{134}, (2020), 109--685. \\
	
	
	\bibitem{Mondal20}
	S. Mondal, M. Biswas, N. Bairagi,
	Local and global dynamics of a fractional-order predator-prey system with habitat complexity and the corresponding discretized fractional-order system.
	\emph{J. Applied Mathematics and Computing.} (2020); 
	DOI: https://doi.org/10.1007/s12190-020-01319-6. \\
	
	
	
	\bibitem{Khan2020}
	A. Khan, J. F. Gomez-Aguilar, T. Abdeljawad,  H. Khan,
	Stability and numerical simulation of a fractional order plant-nectar-pollinator model. \emph{Alexandria Engineering Journal.} \textbf{59}, 1 (2020), 49--59. \\
	
	
	\bibitem{Agarwal13}
	R. P. Agarwal, A. M. A.  El-Sayed, S. M.  Salman,
	Fractional-order Chua's system: discretization, bifurcation and chaos.
	\emph{Advances in Difference Equations.}  \textbf{1}, (2013), 320. \\
	
	\bibitem{Sayed13}
	A. M. A.  El-Sayed, S. M.  Salman,
	On a discretization process of fractional-order Riccati differential equation.
	\emph{J. Fract. Calc. Appl.} \textbf{4}, 2 (2013), 251--259. \\
	
	
	\bibitem{Selvam18} 
	A. G. M. Selvam, D. A. Vianny,
	Behavior of a Discrete Fractional Order SIR Epidemic model.
	\emph{International Journal of Engineering and Technology.} \textbf{7}, 4.10 (2018), 675--680.
	
	\bibitem{Raheem14}
	Z. F. El Raheem,  S. M.  Salman,
	On a discretization process of fractional-order logistic differential equation.
	\emph{Journal of the Egyptian Mathematical Society.} \textbf{22}, 3 (2014), 407--412. 
	
	
	\bibitem{ChattoBairagi01} 
	J. Chattopadhyay, N. Bairagi,
	Pelicans at risk in Salton sea--an eco-epidemiological model.
	\emph{Ecological Modelling.} \textbf{136}, (2001), 103--112.
	
	\bibitem{Elaydi99} 
	S. N. Elaydi,
	\emph{An Introduction to Difference Equations}.
	Springer, New York, NY, USA, (1999).    
	
	
	\bibitem{Bairagi08} 
	N. Bairagi, R. R. Sarkar, J. Chattopadhyay,
	Impacts of incubation delay on the dynamics of an eco-epidemiological system--a theoretical study.
	\emph{Bulletin of mathematical biology.} \textbf{70},  7 (2008), 2017--2038.  
	
	\bibitem{ChattoBairagi08}
	R. K. Upadhyay, N. Bairagi, K. Kundu, J. Chattopadhyay,
	Chaos in eco-epidemiological problem of the Salton Sea and its possible control.
	\emph{Applied Mathematics and Computation.} \textbf{196},  1 (2008), 392--401.
	
	
	
\end{thebibliography}
\end{document}